\documentclass{amsart}
\usepackage{graphicx}
\input xy
\xyoption{all}
\newenvironment{theorem}[2][Theorem]{\begin{trivlist}
\item[\hskip \labelsep {\bfseries #1}\hskip \labelsep {\bfseries #2}]}{\end{trivlist}}

\newtheorem{thm}{Theorem}
\newtheorem{cor}[thm]{Corollary}
\newtheorem{lem}[thm]{Lemma}
\newtheorem{prop}[thm]{Proposition}
\theoremstyle{definition}

\theoremstyle{remark}
\newtheorem{rem}[thm]{Remark}
\theoremstyle{definition}
\newtheorem{ex}[thm]{Examples}
\numberwithin{equation}{section}

\begin{document}

\title[A descent theorem for formal smoothness]{A descent theorem for formal smoothness}%
\author{Javier Majadas}%
\address{Departamento de Matem\'aticas \newline Universidad de Santiago de Compostela \newline E15782 Santiago de Compostela \newline Spain}%
\email{j.majadas@usc.es}%



\begin{abstract}
We give a descent result for formal smoothness having interesting applications: we deduce that quasi-excellence descends along flat local homomorphisms of finite type, we greatly improve Kunz's characterization of regular local rings by means of the Frobenius homomorphisms as well as Andr\'e and Radu relativization of this result, etc. In the second part of the paper, we study a similar question for the complete intersection property instead of formal smoothness, giving also some applications.
\end{abstract}

\maketitle

In this paper we obtain a descent result for formal smoothness (and regularity). This single result (Theorem \ref{13}) has as particular cases a variety of (a priori unconnected) results in commutative algebra. Some of these particular cases greatly improve known important results (the descent of quasi-excellence by finite surjective morphisms \cite {Gr}, Kunz's characterization of regular local rings in positive characteristic \cite {Ku}, \cite {Ro1}, as well as its Andr\'e-Radu relativization \cite {AnCRAS}, \cite {AnManus}, and its analogue for complete intersections \cite {BM}); some other particular cases improve minor results of some papers (\cite{Av-b}, \cite{AHS}, \cite{AHIY}, \cite{Du}, \cite{Du2}, \cite{TY}); some others give already known theorems (\cite{Ro2}) or little results unexplored up to now (Corollary \ref{lastcor}), but the point is that it is the same main result obtained to prove the descent of quasi-excellence the one that gives also those other results. Since this main result is rather technical, we confine ourselves in this introduction to discuss these applications in more detail.\\
\\*
\emph {A. Descent of quasi-excellence.} The descent of quasi-excellence along flat local homomorphisms was studied in \cite{Gr}, where it was shown that quasi-excellence descends along {\it finite} flat local homomorphisms.
We remove here this restriction, proving that this property descends along flat local homomorphisms of finite type. More precisely, Greco in \cite{Gr} deduced this result from the following theorem: if $f:R \rightarrow S$ is a finite local homomorphism of noetherian local rings such that  $Spec (S) \rightarrow Spec (R)$ is surjective, then quasi-excellence of $S$ implies quasi-excellence of $R$. He also obtained an example \cite [Proposition 4.1]{Gr} which shows that this last theorem cannot be extended to the case where $f$ is of finite type instead of finite, even when $R$ and $S$ are local domains of dimension 1 (the problem is concentrated in the local condition, i.e., the geometric regularity of the formal fibers, since the J-2 property always descends and ascends by surjective morphisms of finite type \cite [Corollary 2.4]{Gr}). So the problem for finite type flat homomorphisms was left aside (however he was able to show that at least $R$ contains a non-empty principal open that is quasi-excellent), and he proceeded instead to consider proper surjective morphisms of schemes instead, a problem which was finally solved in \cite {Be}, \cite {Og}.

We reconsider here the problem for homomorphisms of finite type, and we will prove that it holds. Note first that if $\hat{R}$ and $\hat{S}$ are the completions of $R$ and $S$ at their maximal ideals, the surjectivity of $Spec (S) \rightarrow Spec (R)$ is equivalent to the surjectivity of $Spec (\hat{S}) \rightarrow Spec (\hat{R})$ when $f$ is finite (but not necessarily when $f$ is of finite type). Then, a particular case of our main theorem gives:

\begin{theorem}{\ref{excellent}.}
Let $f:R \rightarrow S$ be a local homomorphism essentially of finite type (or more generally $S\otimes_R\hat{R}$ is noetherian) of noetherian local rings such that $Spec (\hat{S}) \rightarrow Spec (\hat{R})$ is surjective. If $S$ is quasi-excellent then $R$ is quasi-excellent. In particular, quasi-excellence descends along flat local homomorphisms (essentially) of finite type.\\
\end{theorem}

\noindent \emph {B. Extension of Kunz's theorem.} In \cite {Ku}, E. Kunz obtained the following result:
\begin{theorem}{(Kunz)}
Let $A$ be a noetherian local ring containing a field of characteristic $p>0$, let $\phi :A \rightarrow A, \phi(a)=a^p$ be the Frobenius homomorphism, and let $^{\phi}\!A$ be the ring $A$ considered as $A$-module via $\phi$. Then $A$ is regular if and only if $\phi: A \rightarrow {^\phi A}$ is flat.\\
\end{theorem}

The relative case was obtained by N. Radu and M. Andr\'e \cite {Ra}, \cite {AnCRAS} and \cite {AnManus}:

\begin{theorem}{(Andr\'e, Radu)} Let $\alpha:A \rightarrow B$ be a homomorphism of noetherian rings containing a field of characteristic $p>0$. Then $\alpha:A \rightarrow B$ is regular if and only if the relative Frobenius homomorphism $\phi_{B|A} :{^{\phi}\!A} \otimes_AB \rightarrow {^{\phi}\!B}$ is flat.\\

\end{theorem}

This theorem contains Kunz's result as a particular case: if $A$ is a noetherian local ring containing a field of characteristic $p>0$, then $A$ contains a perfect field $F$, and applying this theorem to the homomorphism $F \rightarrow A$ we recover Kunz's theorem.

The ``if'' part in the absolute Kunz's theorem was extended in \cite {TY}, \cite {AIM}, \cite {AHIY}, from the particular case of the Frobenius endomorphism to the more general case of a contracting endomorphism, that is, a homomorphism of noetherian local rings $f:(A,\mathfrak{m},k) \rightarrow (A,\mathfrak{m},k)$ such that there exists some $i>0$ with $f^i(\mathfrak{m})\subset \mathfrak{m}^2$.\\

We show how our main result has also as a special case  a relative version of these results for contracting homomorphisms, that is, we obtain a result (Theorem \ref{theorem}) generalizing at once both the relative result of Andr\'e and Radu and the above extensions for contracting homomorphisms of the absolute result by Kunz.\\
\\*
\emph {C. Rodicio's theorem.}
\begin{theorem}{(Rodicio)}
Let $u:A \rightarrow B$ be a flat homomorphism of noetherian rings such that $B\otimes_AB$ is noetherian. If the flat (Tor) dimension of $B$ over $B\otimes_AB$ is finite, then $u$ is regular.
\end{theorem}

This is the main result in \cite {Ro2}. Previous results on this topic were obtained by Auslander, Eilenberg, Harada, Hochschild, Nakayama, Rosenberg, Zelinsky and Rodicio himself among others. Our main result also has this theorem as a special case (Theorem \ref{Rodicio}). However, it should be noted that Rodicio's proof is valid even when $B\otimes_AB$ is not noetherian.\\
\\*
\emph {D. Decomposition of formal smoothness.} Another particular case of our main result is the following one (Theorem \ref{decomposition}). Let $A \xrightarrow{u} B \xrightarrow{v} C$ be local homomorphisms of noetherian local rings such that the flat dimension of $C$ over $B$ is finite; if $vu$ is formally smooth then $u$ is formally smooth.

When $C$ is a flat $B$-module (i.e., the flat dimension is zero), the result is well known (it is an easy consequence of \cite [0$_{IV}$19.7.1, 0$_{IV}$22.5.8]{EGAIV1}). For finite flat dimension it follows from \cite [Theorem 4.7]{Av-b}.\\
\\*
\emph {E.} There are also a number of particular cases of our main result that, though important, they are well known and also easy to prove by using Andr\'e-Quillen homology. We only point out a couple of them in order to stress the ubiquity of our main descent theorem (Examples \ref{other}).\\

In the second section, we study similar questions for the complete intersection property instead of formal smoothness. For it, we need first some analogues of the main results of \cite {Av-b} for finite complete intersection dimension instead of finite flat dimension (Theorem \ref{TheoremA} - Remark \ref{remC}). We now give a more reduced list of applications, since many of the ones given for formal smoothness are easy to translate to the complete intersection case.\\

We will use some facts of Andr\'e-Quillen homology modules $H_n(A,B,M)$. For convenience of the reader we list them in an Appendix with precise references.\\

\section{Formal smoothness}

\begin{thm}\label{13}
Let
$$
\xymatrix{  (A,\mathfrak{m},k) \ar[r]^u \ar[d]^f & (B,\mathfrak{n},l) \ar[d]^g  \\
(\tilde{A},\mathfrak{\tilde{m}},\tilde{k}) \ar[r]^{\tilde{u}}  & (\tilde{B},\mathfrak{\tilde{n}},\tilde{l}) }
$$
be a commutative square of local homomorphisms of noetherian local rings such that\\*
(i) $Tor_i^A(\tilde{A},B)=0$ for all $i>0$.\\*
(ii) The homomorphism $H_1(A,B,\tilde{l}) \rightarrow H_1(\tilde{A},\tilde{B},\tilde{l})$ vanishes.\\*
(iii) If $\mathfrak{p}$ is the contraction in $\tilde{A}\otimes_AB$ of the maximal ideal $\mathfrak{\tilde{n}}$ of $\tilde{B}$, then $(\tilde{A}\otimes_AB)_{\mathfrak{p}}$ is a noetherian ring.\\*
(iv) $fd_{\tilde{A}\otimes_AB}(\tilde{B}) <\infty$.\\

Then $u:A \rightarrow B$ is formally smooth.
\end{thm}
\begin{proof}
\noindent  The commutative diagram of ring homomorphisms
$$
\xymatrix{  \tilde{A} \ar[r] \ar[d]^{=} & \tilde{A}\otimes_AB \ar[r] \ar[d]^{=} & \tilde{B} \ar[d] \\
\tilde{A} \ar[r]  & \tilde{A}\otimes_AB \ar[r]  & \tilde{l}  }
$$
gives, taking Jacobi-Zariski sequences, a commutative square

$$
\xymatrix{  H_2(\tilde{A}\otimes_AB,\tilde{B},\tilde{l}) \ar[r]^\delta \ar[d]^{\gamma} & H_1(\tilde{A},\tilde{A}\otimes_AB,\tilde{l}) \ar[d]^{=}  \\
H_2(\tilde{A}\otimes_AB,\tilde{l},\tilde{l}) \ar[r]  & H_1(\tilde{A},\tilde{A}\otimes_AB,\tilde{l}). }
$$
Associated to the ring homomorphisms $\tilde{A}\otimes_AB \rightarrow \tilde{B} \rightarrow \tilde{l}$ we have an exact sequence
$$H_2(\tilde{A}\otimes_AB,\tilde{B},\tilde{l}) \xrightarrow{\gamma} H_2(\tilde{A}\otimes_AB,\tilde{l},\tilde{l})= H_2((\tilde{A}\otimes_AB)_{\mathfrak{p}},\tilde{l},\tilde{l}) \xrightarrow{\beta} H_2(\tilde{B},\tilde{l},\tilde{l}).$$
where $\mathfrak{p}$ is the contraction in $\tilde{A}\otimes_AB$ of the maximal ideal of $\tilde{B}$. Avramov's result (Appendix 11) says that $\beta$ is injective. Thus $\gamma=0$.

On the other hand, since $Tor_i^A(\tilde{A},B)=0$ for all $i>0$, the canonical homomorphism $H_1(A,B,\tilde{l}) \rightarrow H_1(\tilde{A},\tilde{A}\otimes_AB,\tilde{l})$ is an isomorphism, so the commutative triangle
$$
\xymatrix{  H_1(A,B,\tilde{l}) \ar[dr]^0 \ar[rr]^{\simeq} && H_1(\tilde{A},\tilde{A}\otimes_AB,\tilde{l}) \ar[dl]^{\alpha} \\
& H_1(\tilde{A},\tilde{B},\tilde{l}) & }
$$
shows that $\alpha =0$. From the ring homomorphisms $\tilde{A} \rightarrow \tilde{A}\otimes_AB \rightarrow \tilde{B}$, we obtain an exact sequence
$$H_2(\tilde{A}\otimes_AB,\tilde{B},\tilde{l}) \xrightarrow{\delta} H_1(\tilde{A},\tilde{A}\otimes_AB,\tilde{l}) \xrightarrow{\alpha} H_1(\tilde{A},\tilde{B},\tilde{l})$$
and then $\delta$ is surjective.

Thus, from the above commutative square we deduce that $H_1(\tilde{A},\tilde{A}\otimes_AB,\tilde{l})=0$, and so $H_1(A,B,\tilde{l})=0$. This implies that $u:A \rightarrow B$ is formally smooth.

\end{proof}

\begin{rem} (i) In the particular case where $f$ and $g$ are the Frobenius endomorphisms, this result was proved in \cite{Du2}.\\*
(ii) From the proof of Theorem \ref{13}, it is also clear that we can also weaken the condition $fd_{\tilde{A}\otimes_AB}(\tilde{B}) <\infty$ as follows: there exists a local homomorphism of local noetherian rings $\tilde{B} \rightarrow \tilde{C}$ such that the composition map $fd_{\tilde{A}\otimes_AB}(\tilde{C}) <\infty$.
\end{rem}

In order to apply this result, it will be convenient to state some corollaries first:

\begin{cor}\label{last}
Let
$$
\xymatrix{  (A,\mathfrak{m},k) \ar[r]^u \ar[d]^f & (B,\mathfrak{n},l) \ar[d]^g  \\
(\tilde{A},\mathfrak{\tilde{m}},\tilde{k}) \ar[r]^{\tilde{u}}  & (\tilde{B},\mathfrak{\tilde{n}},\tilde{l}) }
$$
be a commutative square of local homomorphisms of noetherian local rings such that\\*
(i) $Tor_i^A(\tilde{A},B)=0$ for all $i>0$.\\*
(ii) The homomorphism $H_1(A,B,\tilde{l}) \rightarrow H_1(\tilde{A},\tilde{B},\tilde{l})$ is the zero map.\\*
(iii) The $\tilde{A}\otimes_AB$-module $\tilde{B}$ is flat.\\

Then $u:A \rightarrow B$ is formally smooth.
\end{cor}

\begin{proof}

Condition (iii) of Theorem \ref{13} holds by faithfully flat descent.

\end{proof}

\begin{cor}\label{regnaive}
Let
$$
\xymatrix{  A \ar[r]^u \ar[d]^f & B \ar[d]^g  \\
\tilde{A} \ar[r]^{\tilde{u}}  & \tilde{B} }
$$
be a commutative square of homomorphisms of noetherian rings verifying\\*
(i) $Tor_i^A(\tilde{A},B)=0$ for all $i>0$. \\*
(ii) For each $\mathfrak{q} \in Spec(B)$ there exists $\mathfrak{\tilde{q}} \in Spec(\tilde{B})$ with $g^{-1}(\mathfrak{\tilde{q}}) = \mathfrak{q}$ such that the map
$$H_1(A,B,k(\mathfrak{\tilde{q}})) \rightarrow H_1(\tilde{A},\tilde{B},k(\mathfrak{\tilde{q}}))$$
vanishes, where $k(\mathfrak{\tilde{q}})$ is the residue field of $\tilde{B}_{\mathfrak{\tilde{q}}}$.\\*
(iii) The $\tilde{A}\otimes_AB$-module $\tilde{B}$ is flat.\\

Then the homomorphism $u:A \rightarrow B$ is regular, that is, flat with geometrically regular fibers.
\end{cor}

\begin{proof}
By \cite [0$_{IV}$19.7.1, 22.5.8]{EGAIV1}, we have to prove that the local homomorphism $A_{\mathfrak{p}} \rightarrow B_{\mathfrak{q}}$ is formally smooth for each $\mathfrak{q} \in Spec(B)$, where $\mathfrak{p} = u^{-1}(\mathfrak{q})$.
Let $\mathfrak{\tilde{q}} \in Spec(\tilde{B})$ such that $g^{-1}(\mathfrak{\tilde{q}}) = \mathfrak{q}$ as in condition (ii), and let $\mathfrak{\tilde{p}} = \tilde{u}^{-1}(\mathfrak{\tilde{q}})$. We have\\*
(i) $Tor_i^{A_{\mathfrak{p}}}(\tilde{A}_{\mathfrak{\tilde{p}}},B_{\mathfrak{q}})=0$ for all $i>0$.\\*
(ii) The homomorphism
$$H_1(A_{\mathfrak{p}},B_{\mathfrak{q}},k(\mathfrak{\tilde{q}})) = H_1(A,B,k(\mathfrak{\tilde{q}})) \rightarrow H_1(\tilde{A},\tilde{B},k(\mathfrak{\tilde{q}})) = H_1(\tilde{A}_{\mathfrak{\tilde{p}}},\tilde{B}_{\mathfrak{\tilde{q}}},k(\mathfrak{\tilde{q}}))$$
is zero.\\*
(iii) The $\tilde{A}_{\mathfrak{\tilde{p}}}\otimes_{A_{\mathfrak{p}}}B_{\mathfrak{q}}$-module $\tilde{B}_{\mathfrak{\tilde{q}}}$ is flat.\\

Therefore we can apply Corollary \ref{last} to the commutative square
$$
\xymatrix{  A_{\mathfrak{p}} \ar[r] \ar[d] & B_{\mathfrak{q}} \ar[d]  \\
\tilde{A}_{\mathfrak{\tilde{p}}} \ar[r]  & \tilde{B}_{\mathfrak{\tilde{q}}}}
$$
and we obtain that $A_{\mathfrak{p}} \rightarrow B_{\mathfrak{q}}$ is formally smooth.

\end{proof}

In the applications, we will use frequently the following weaker version of the theorem:

\begin{cor}\label{14}
Let
$$
\xymatrix{  (A,\mathfrak{m},k) \ar[r]^u \ar[d]^f & (B,\mathfrak{n},l) \ar[d]^g  \\
(\tilde{A},\mathfrak{\tilde{m}},\tilde{k}) \ar[r]^{\tilde{u}}  & (\tilde{B},\mathfrak{\tilde{n}},\tilde{l}) }
$$
be a commutative square of local homomorphisms of noetherian local rings such that\\*
(i) $Tor_i^A(\tilde{A},B)=0$ for all $i>0$.\\*
(ii) The homomorphism $\tilde{u}:\tilde{A} \rightarrow \tilde{B}$ is formally smooth.\\*
(iii) The ring $\tilde{A}\otimes_AB$ is noetherian.\\*
(iv) $fd_{\tilde{A}\otimes_AB}(\tilde{B}) <\infty$.\\

Then $u:A \rightarrow B$ is formally smooth.
\end{cor}

\begin{proof}
Condition (ii) means $H_1(\tilde{A},\tilde{B},\tilde{l})=0$.
\end{proof}

\begin{rem}\label{factorization}
Another particular case of Theorem \ref{13}, a little more general than Corollary \ref{14}, can be useful. It consists in substituting in Corollary \ref{14} condition (ii) with:\\*
(ii') The square
$$
\xymatrix{  (A,\mathfrak{m},k) \ar[r]^u \ar[d]^f & (B,\mathfrak{n},l) \ar[d]^g  \\
(\tilde{A},\mathfrak{\tilde{m}},\tilde{k}) \ar[r]^{\tilde{u}}  & (\tilde{B},\mathfrak{\tilde{n}},\tilde{l}) }
$$
admits a factorization
$$
\xymatrix{  (A,\mathfrak{m},k) \ar[r]^u \ar[d] & (B,\mathfrak{n},l) \ar[d]  \\
 (\bar{A},\mathfrak{\bar{m}},\bar{k}) \ar[r]^{\bar{u}} \ar[d] & (\bar{B},\mathfrak{\bar{n}},\bar{l}) \ar[d]  \\
(\tilde{A},\mathfrak{\tilde{m}},\tilde{k}) \ar[r]^{\tilde{u}}  & (\tilde{B},\mathfrak{\tilde{n}},\tilde{l}) }
$$
with $\bar{u}$ formally smooth.
\end{rem}

\noindent \emph {Applications}\\
\\*
\emph {A. Extension of Greco's theorem.}

Let $u:R \rightarrow S$ be a local homomorphism of finite type of noetherian local rings. If $R$ is quasi-excellent then $S$ is quasi-excellent. We are concerned here with the reciprocal. In \cite{Gr}, Greco shows that if $u$ is finite and $Spec (S) \rightarrow Spec (R)$ is surjective then $S$ quasi-excellent implies that $R$ is quasi-excellent. He also proves that ``$u$ finite'' cannot be replaced by ``$u$ of finite type'', even when $R$ and $S$ are local domains of dimension 1 \cite [Proposition 4.1]{Gr}.

Instead of the surjectivity of $Spec (S) \rightarrow Spec (R)$, we consider here the surjectivity of $Spec (\hat{S}) \rightarrow Spec (\hat{R})$. When $u$ is finite, both conditions are equivalent (since in this case $Spec (\hat{S}) \rightarrow Spec (\hat{R})$ is obtained from $Spec (S) \rightarrow Spec (R)$ by base change and surjectivity of a morphism of schemas is a property stable by base change). However, for a homomorphism of finite type $u$, the latter condition is in general strictly stronger. Then we can prove:

\begin{thm}\label{excellent}
Let $u:R \rightarrow S$ be a local homomorphism essentially of finite type of noetherian local rings such that $Spec (\hat{S}) \rightarrow Spec (\hat{R})$ is surjective. If $S$ is quasi-excellent then $R$ is quasi-excellent.
\end{thm}

\begin{proof}
A noetherian local ring $A$ is quasi-excellent if and only if the homomorphism $A \rightarrow \hat{A}$ is regular. So we apply Corollary \ref{regnaive} to the square
$$
\xymatrix{  R \ar[r] \ar[d]_u & \hat{R} \ar[d]  \\
S \ar[r]  & \hat{S} }
$$
and use e.g. \cite [7.9.3.1]{EGAIV2} to see that $S\otimes_R \hat{R} \rightarrow \hat{S}$ is flat.

\end{proof}

\begin{cor}\label{flatexcellent}
Let $R$ be a local ring and $u:R \rightarrow S$ be a local flat homomorphism essentially of finite type. If $S$ is quasi-excellent then $R$ is quasi-excellent.

\end{cor}
\begin{proof}
Since $S$ is noetherian, so is $R$. The local homomorphism $\hat{R} \rightarrow \hat{S}$ is flat by the local flatness criterion and so $Spec (\hat{S}) \rightarrow Spec (\hat{R})$ is surjective.
\end{proof}

\noindent \emph {B. Extension of Kunz's theorem.}

We will see now that another special case of Theorem \ref{13} gives us a relative version of Kunz's result for arbitrary contracting homomorphisms (Theorem \ref{theorem}; see the Introduction for comments). In order to guarantee that these homomorphisms satisfy the hypotheses of our theorem, first we will need to prove some facts. We start by noticing that \cite [Lemma 1]{Du} holds also for contracting homomorphisms.

\begin{lem}\label{flat}
Let
$$
\xymatrix{  (A,\mathfrak{m},k) \ar[r]^u \ar[d]^f & (B,\mathfrak{n},l) \ar[d]^g  \\
(A,\mathfrak{m},k) \ar[r]^{u}  & (B,\mathfrak{n},l) }
$$
be a commutative square of local homomorphisms of noetherian local rings such that $f(\mathfrak{m}) \subset \mathfrak{m}^2$. If the homomorphism $\omega:A\otimes_AB \rightarrow B$, $\omega(a\otimes b)=u(a)g(b)$, is flat then $u$ is flat.
\end{lem}

\begin{proof}
By the local flatness criterion it suffices to show that
$$ u_n := A/\mathfrak{m}^n\otimes_Au : A/\mathfrak{m}^n \rightarrow B/\mathfrak{m}^nB$$
is flat for all $n \geq 1$.

We recall from the introduction that ${^{f}\!A}$ is the ring $A$ considered as $A$-module via $f: A \rightarrow A$ and similarly ${^{g}\!B}$, etc. The homomorphism $u={^{f}\!u}$ factorizes as
$${^{f}\!A} \xrightarrow{{^{f}\!A}\otimes_Au} {^{f}\!A}\otimes_AB \xrightarrow{w} {^{g}\!B},$$
and so $u_n = {^{f}\!u_n} = {^{f}\!(A/\mathfrak{m}^n)}\otimes_{^{f}\!A}{^{f}\!u}$ factorizes as
$${^{f}\!(A/\mathfrak{m}^n)} \xrightarrow{{^{f}\!(A/\mathfrak{m}^n)}\otimes_{A/\mathfrak{m}^n}u_n} {^{f}\!(A/\mathfrak{m}^n)} \otimes_{A/\mathfrak{m}^n}B/\mathfrak{m}^nB \xrightarrow{\omega_n} {^{g}\!(B/\mathfrak{m}^nB)}$$
where ${\omega_n}$ is defined in the obvious way as the natural map from the pushout of ${^{f}\!(A/\mathfrak{m}^n)}$ and $B/\mathfrak{m}^nB$ over ${A/\mathfrak{m}^n}$.

But for $n \geq 2$ we have $${^{f}\!(A/\mathfrak{m}^n)}\otimes_{A/\mathfrak{m}^n}u_n = {^{f}\!(A/\mathfrak{m}^n)}\otimes_{A/\mathfrak{m}^{n-1}}(A/\mathfrak{m}^{n-1}\otimes_{A/\mathfrak{m}^n}u_n) = {^{f}\!(A/\mathfrak{m}^n)}\otimes_{A/\mathfrak{m}^{n-1}}u_{n-1}$$
(we have used that $f(\mathfrak{m}) \subset \mathfrak{m}^2$). The homomorphism $u_1$ is flat  since $A/\mathfrak{m}$ is a field. By induction, if $u_{n-1}$ is flat, then ${^{f}\!(A/\mathfrak{m}^n)}\otimes_{A/\mathfrak{m}^n}u_n$ is flat and thus $u_n=\omega_n \circ ({^{f}\!(A/\mathfrak{m}^n)}\otimes_{A/\mathfrak{m}^n}u_n)$ is flat, since $\omega_n = {^{f}\!(A}/\mathfrak{m}^n)\otimes_{{^{f}\!A}}\omega$ is flat for all $n \geq 1$ by base change.

\end{proof}

Now a version of \cite [10.11]{An1974} for contracting endomorphisms.
\begin{lem}\label{Tor}
Let $A$ be a noetherian ring, $I$ an ideal of $A$, $M$ an $A$-module of finite type, $f:A \rightarrow A$ a ring homomorphism such that $f(I) \subset I^2$, $g: M\rightarrow M$ an $f$-homomorphism ($g(am) = f(a)g(m)$ for $a\in A, m \in M$), and $n>0$ an integer. Then there exists an $s>0$ such that the map induced by $f^s$ and $g^s$ (in all three variables)
$$ Tor^A_n(M,A/I) \rightarrow Tor^A_n(M,A/I)$$
is zero.
\end{lem}

\begin{proof}
Let
$$ ... \rightarrow F_2 \xrightarrow{d_2} F_1 \xrightarrow{d_1} F_0 \rightarrow M \rightarrow 0$$
be a resolution of $M$ with $F_i$ a projective $A$-module of finite type for each $i$. Since
$$Tor^A_n(M,A/I) = H_n(F_*/IF_*)$$
is the homology of the complex
$$F_{n+1}/IF_{n+1} \xrightarrow{\delta_{n+1}} F_{n}/IF_{n} \xrightarrow{\delta_{n}} F_{n-1}/IF_{n-1}$$
we have
$$Tor^A_n(M,A/I) = Ker (\delta_{n}) / Im (\delta_{n+1}) = Ker \left(\frac{F_{n}/IF_{n}}{Im (\delta_{n+1})} \rightarrow F_{n-1}/IF_{n-1}\right).$$

Applying $-\otimes_AA/I$ to the exact sequence
$$ 0 \rightarrow Im (d_{n+1}) \rightarrow F_n \rightarrow Im (d_n) \rightarrow 0$$
we get
$$\frac{F_{n}/IF_{n}}{Im (\delta_{n+1})} = \frac{Im(d_n)}{I\cdot Im(d_n)}.$$

Thus
$$Tor^A_n(M,A/I) = Ker \left(\frac{Im(d_n)}{I\cdot Im(d_n)} \rightarrow F_{n-1}/IF_{n-1}\right) = \frac{Im(d_n)\cap  IF_{n-1}}{I\cdot Im(d_n)}.$$

By the Artin-Rees lemma, there exists a positive integer $t$ such that $Im(d_n)\cap  I^tF_{n-1} \subset I\cdot Im(d_n)$, so choosing $s$ such that $f^s(I) \subset I^t$ we see that $f^s$ and $g^s$ induce the zero map on
$$Tor^A_n(M,A/I) = \frac{Im(d_n)\cap  IF_{n-1}}{I\cdot Im(d_n)}.$$

\end{proof}

\begin{prop}\label{vanish}
Let $(A,\mathfrak{m},k)$ be a noetherian local ring and $f:A \rightarrow A$ a ring homomorphism such that $f(\mathfrak{m}) \subset \mathfrak{m}^2$. For each integer $n\geq 0$ there exists an integer $s>0$ such that $f^s$ induces the zero map of functors
$$H_n(A,k,-) \rightarrow H_n(A,k,-).$$
\end{prop}

\begin{proof}
We mean that for any $k$-module $M$, the map induced by $f^s$ on the left two variables $H_n(A,k,M) \rightarrow H_n(A,k,M)$ is zero, where in the left $H_n(A,k,M)$, $M$ is considered as $k$-module by restriction of scalars via the map $k \rightarrow k$ induced by $f^s$.

We have $H_0(A,k,k)=0$, so we can assume $n>0$. By Lemma \ref{Tor},  for $i=1, ... n$, we can choose integers $s_i$ such that $f^{s_i}$ induces the zero map  $Tor^A_i(k,k) \rightarrow Tor^A_i(k,k)$. We will use \cite [10.12]{An1974} with $A^i=A$, $B^i=C^i=k$ for all $i$, and the homomorphisms $A^{i-1} \rightarrow A^i$, $B^{i-1} \rightarrow B^i$, $C^{i-1} \rightarrow C^i$ are the ones induces by $f^{s_i}$. We obtain that for $s \geq s_1+...+s_n$, the map $H_n(A,k,-) \rightarrow H_n(A,k,-)$ induced by $f^s$ is zero.
\end{proof}

\begin{cor}\label{corollary4}
Let
$$
\xymatrix{  (A,\mathfrak{m},k) \ar[r]^u \ar[d]^f & (B,\mathfrak{n},l) \ar[d]^g  \\
(A,\mathfrak{m},k) \ar[r]^u  & (B,\mathfrak{n},l) }
$$
be a commutative square of local homomorphisms of noetherian local rings such that $f(\mathfrak{m}) \subset \mathfrak{m}^2$, $g(\mathfrak{n}) \subset \mathfrak{n}^2$. Then for each integer $n\geq 0$ there exists an integer $s>0$ such that $(f^s,g^s)$ induces the zero homomorphism
$$H_n(A,B,M) \xrightarrow{\alpha_s} H_n(A,B,M)$$
for any $l$-module $M$.
\end{cor}
\begin{proof}
We have a commutative diagram with exact rows

$$
\xymatrix{  H_{n+1}(B,l,-) \ar[r]^\delta \ar[d]^{\beta_i} & H_n(A,B,-) \ar[r]^{\varepsilon} \ar[d]^{\alpha_i} & H_n(A,l,-) \ar[d]^{\gamma_i} \\
H_{n+1}(B,l,-) \ar[r]^\delta  & H_n(A,B,-) \ar[r]^{\varepsilon}  & H_n(A,l,-)  }
$$
where the vertical maps are induced by $(f^i,g^i)$. If $\beta_i=0$, $\gamma_i=0$, then $Im(\alpha_i) \subset Ker(\varepsilon) = Im(\delta)$ and $Im(\delta) \subset Ker(\alpha_i)$, so $\alpha_{2i}=\alpha_i^2=0$. Therefore, taking $i$ sufficiently large, Proposition \ref{vanish} gives the result.
\end{proof}

Now, we can see that the required extension of Kunz's theorem is also a particular case of our Theorem \ref{13}. Since we are placed in the local case, we do not need the flatness of $\omega$, but only the flatness of its completion at the desired prime.

\begin{thm}\label{theorem}
Let
$$
\xymatrix{  (A,\mathfrak{m},k) \ar[r]^u \ar[d]^f & (B,\mathfrak{n},l) \ar[d]^g  \\
(A,\mathfrak{m},k) \ar[r]^u  & (B,\mathfrak{n},l) }
$$
be a commutative square of local homomorphisms of noetherian local rings such that there exists some $i$ with $f^i(\mathfrak{m}) \subset \mathfrak{m}^2$, $g^i(\mathfrak{n}) \subset \mathfrak{n}^2$.

Denote by $\hat{A}$ and $\hat{B}$ the completions of $A$ and $B$ at their maximal ideals. If there exists some $j$ such that the homomorphism $\hat{\omega}_j:{^{f_j}\!\hat{A}}\otimes_{\hat{A}} \hat{B} \rightarrow {^{g_j}\!\hat{B}}$ is flat, then $u: A \rightarrow B$ is formally smooth.
\end{thm}

\begin{proof}
By \cite [0$_{IV}$19.3.6]{EGAIV1}, the homomorphism $u: A \rightarrow B$ is formally smooth if (and only if) $u: \hat{A} \rightarrow \hat{B}$ is. So we can assume that $\omega_j:{^{f_j}\!A}\otimes_{A} B \rightarrow {^{g_j}\!B}$ is flat. By Corollary \ref{corollary4} there exists an integer $s$ such that $(f^{is},g^{is})$ induces the zero homomorphism $H_1(A,B,M) \rightarrow H_1(A,B,M)$ for any $l$-module $M$. On the other hand, if $\omega_j$ is flat, then $\omega_{2j}$ is also flat, since it coincides with the composition
$${^{f^{2j}}\!A}\otimes_AB \xrightarrow{{^{f^{2j}}\!A}\otimes_{^{f^{j}}\!A}\omega_j} {^{f^{2j}}\!A}{\otimes_{^{f^{j}}\!A}} ^{g^{j}}\!B \xrightarrow{\omega_j} {^{g^{2j}}\!B}$$
of two flat homomorphisms.
So, replacing $(f,g)$ by a suitable power $(f^{t},g^{t})$ we can assume $f(\mathfrak{m}) \subset \mathfrak{m}^2$, $\omega : {^{f}\!A}\otimes_AB \rightarrow {^g\!B}$ flat, and the homomorphism $(f,g)$ induces the zero map $H_1(A,B,M) \rightarrow H_1(A,B,M)$.

By Lemma \ref{flat}, $u$ is flat. Then the result follows from Corollary \ref{last}.
\end{proof}

\begin{rem}\label{rem.13.5} As a particular case of Theorem \ref{theorem}, we have obviously the original Kunz's theorem (even for finite flat dimension instead of flatness if we use Theorem \ref{13}): let $A$ be a noetherian local ring containing a field of characteristic $p>0$, let $\phi_A :A \rightarrow A$ be the Frobenius homomorphism. If $fd_A({^{\phi_A}A})<\infty$ then $A$ is regular. This result follows from Theorem \ref{13} and Corollary \ref{corollary4} applied to the commutative square
$$
\xymatrix{  k \ar[r] \ar[d]^{\phi_k^s} & A \ar[d]^{\phi_A^s}  \\
k \ar[r] & A }
$$
for some $s$ sufficiently large, where $k$ is a perfect field contained in $A$. When $fd_A({^{\phi_A}A})=0$ this is Kunz's result \cite{Ku}, and in general it was obtained by Rodicio \cite{Ro1}. If we consider an arbitrary contracting endomorphism $\phi$ instead of $\phi_A$, taking as $k$ a perfect field contained in the subfield of elements of $A$ fixed by $\phi$, we obtain a similar result, which by different methods was obtained previously in \cite {TY}, \cite {AHIY}.

Finally, in the general setting, we may ask if our Theorem \ref{13} gives also a relativization of \cite [Proposition 2]{RCIG}. That is, for a local homomorphism $f:A \rightarrow B$, we ask if \cite [Proposition 2]{RCIG} is a consequence of our Theorem \ref{13} applied to
$$
\xymatrix{  \mathbb{Z} \ar[r] \ar @{=}[d] & A \ar[d]^{f}  \\
\mathbb{Z} \ar[r] & B }
$$
In general, the answer is no, since the homology vanishing hypothesis required for $f:A \rightarrow B$ in \cite {RCIG} is weaker than ours. However, when $A$ contains a field, taking as before a perfect subfield $k$ instead of $\mathbb{Z}$, we can prove that the answer is affirmative.

\end{rem}

\noindent \emph {C. Rodicio's theorem.}

\begin{thm}\label{Rodicio}
Let $u:A \rightarrow B$ be a flat homomorphism of noetherian rings such that $B\otimes_AB$ is noetherian. If $fd_{B\otimes_AB}(B)<\infty$ (via the multiplication map $\mu :B\otimes_AB \rightarrow B$) then $u$ is regular.
\end{thm}

\begin{proof}
This is the particular case $\tilde{A}=\tilde{B}=B$. More precisely, apply Corollary \ref{14} to the commutative square
$$
\xymatrix{  A_{\mathfrak{p}} \ar[r] \ar[d] & B_{\mathfrak{q}} \ar @{=}[d]  \\
B_{\mathfrak{q}} \ar @{=}[r]  & B{_\mathfrak{q}} }
$$
\end{proof}

In this generality, this result was first proved by Rodicio \cite{Ro2} (even when $B\otimes_AB$ is not noetherian), and subsequently improved (in the noetherian case) by several authors.

\begin{rem}\label{augmented}
We can see this result as a particular case of Corollary \ref{14} also in a different way. Instead of the above square, we use
$$
\xymatrix{  B_{\mathfrak{q}} \ar[r] \ar @{=}[d] & (B \otimes_AB)_{\mathfrak{c}} \ar[d]  \\
B_{\mathfrak{q}} \ar @{=}[r]  & B{_\mathfrak{q}} }
$$
where $\mathfrak{c} = \mu ^{-1}(\mathfrak{q})$. We deduce that $(B \otimes_AB)_{\mathfrak{c}}$ is formally smooth over $B{_\mathfrak{q}}$ and so
$$H_2((B \otimes_AB)_{\mathfrak{c}},B{_\mathfrak{q}},B{_\mathfrak{q}}/\mathfrak{q}B{_\mathfrak{q}}) = H_1(B{_\mathfrak{q}},(B \otimes_AB)_{\mathfrak{c}},B{_\mathfrak{q}}/\mathfrak{q}B{_\mathfrak{q}})=0,$$
that is, $ker(\mu)$ is locally generated by a regular sequence.

The value of this second approach is that it can be applied to an arbitrary noetherian supplemented $B$-algebra $S$ instead of $B \otimes_AB$.\\

\end{rem}

\noindent \emph {D. Decomposition of formal smoothness.}

\begin{thm}\label{decomposition}
Let $A \xrightarrow{u} B \xrightarrow{v} C$ be local homomorphisms of noetherian local rings such that $vu$ is formally smooth and $fd_B(C)<\infty$. Then $u$ is formally smooth.
\end{thm}
\begin{proof}
This result is Corollary \ref{14} for $\tilde{A}=A$, $\tilde{B}=C$, that is, for the square
$$
\xymatrix{  A \ar[r] \ar @{=}[d] & B \ar[d]  \\
A \ar[r]  & C }
$$
\end{proof}

\noindent \emph {E. Other results.}

There are also other (well known) interesting particular cases of our main result, but they also can be easily proved with the help of Andr\'e-Quillen homology. So we only state here a couple of them:

\begin{ex}\label{other}
(a) 3. Let $(A,\mathfrak{m},k)\rightarrow(B,\mathfrak{n},l)$ be a flat local homomorphism of noetherian local rings. If $k \rightarrow B\otimes_Ak$ is formally smooth then $A \rightarrow B$ is formally smooth.

It follows from Corollary \ref{14} applied to

$$
\xymatrix{  A \ar[r] \ar[d] & B \ar[d]  \\
k \ar[r]  & B\otimes_Ak }
$$

This is part of \cite [0$_{IV}$19.7.1]{EGAIV1} and it is also an immediate consequence of Appendix 4, 10.\\

(b) Let $u:(A,\mathfrak{m},k)\rightarrow(B,\mathfrak{n},l)$ be a flat local homomorphism of noetherian local rings. If $l|k$ is separable and $B\otimes_Ak$ is regular, then $u$ is formally smooth.

We use the same square as in (a) with $B\otimes_Ak$ replaced by $l$. It can be also deduced from the Jacobi-Zariski exact sequence associated to $k \rightarrow B\otimes_Ak \rightarrow l$ and Appendix 4, 9, 10.

\end{ex}

\section{Complete intersection}

We will now prove similar results for complete intersection instead of formal smoothness (Theorem \ref{cidim} and its corollaries). We will need first a different version (Theorems \ref{TheoremA} and \ref{TheoremB}) of two of the main theorems in \cite{Av-b}. For convenience of the reader, we include here the details, though they are easy consequences of Avramov's results.

Avramov \cite {Av} introduces the virtual projective dimension, vpd, as follows (assuming for simplicity that the residue field of $A$ is infinite). An $A$-module $M$ of finite type is said to be of \textit{finite virtual projective dimension} if there exists  a surjective homomorphism of noetherian local rings $Q \rightarrow \hat{A}$ with kernel generated by a regular sequence such that the projective dimension pd$_Q(\hat{A} \otimes_A M)$ is finite. The local ring $A$ is complete intersection if and only if any module of finite type has finite virtual projective dimension if and only if its residue field has finite virtual projective dimension.

Subsequently a modification of this concept, the \textit{complete intersection dimension}, CI-dim, was introduced in \cite {AGP}. In its definition, instead of the completion homomorphism $A \rightarrow \hat{A}$, an arbitrary flat local homomorphism of noetherian local rings $A \rightarrow A'$ is allowed. Complete intersection dimension shares many properties with virtual projective dimension, in particular the above equivalences. It also behaves well with respect to localization, but this advantage can also be considered as a symptom of some difficulties since complete intersection property of homomorphisms does not localize in general. If $f:(A,\mathfrak{m},k) \rightarrow (B,\mathfrak{n},l)$ is a complete intersection homomorphism at the maximal ideal of $B$ and $A$ has complete intersection formal fibers (for instance if $A$ is a quotient of a local complete intersection ring; see \cite {So}), then $f$ is complete intersection at all primes \cite {Mar}, \cite {Tab}, \cite [5.12]{Av-b}. But in general, this property does not localize as we can see taking as $f$ the completion homomorphism of a local ring whose formal fibers are not complete intersection.

For this reason, we consider here a different definition of complete intersection dimension introduced in \cite {So}, \textit{cidim}, that localizes when the ring has complete intersection formal fibers. Moreover, we have vpd $<\infty$ $\Rightarrow$ CI*-dim $<\infty$ $\Rightarrow$ cidim $<\infty$ $\Rightarrow$ CI-dim $<\infty$, where CI*-dim is the upper complete intersection dimension introduced in \cite {T}.

Let $f:(A,\mathfrak{m},k)\rightarrow(R,\mathfrak{n},l)$ be a local homomorphism  of noetherian local rings. We say that $f$ is \textit{weakly regular} if it is flat and the closed fiber $R\otimes_Ak$ is a regular ring. We say that $f$ is \textit{flat complete intersection at $\mathfrak{n}$} if it is flat and the closed fiber $R\otimes_Ak$ is a complete intersection ring. Since we will consider complete intersection homomorphisms always at the maximal ideal, we will simply say flat complete intersection homomorphism.

If $f:(A,\mathfrak{m},k)\rightarrow(B,\mathfrak{n},l)$ is a local homomorphism of noetherian local rings, a \textit{regular factorization} (resp. \textit{complete intersection factorization}) of $f$ is a factorization $A\xrightarrow{i} R \xrightarrow{p} B$ of $f$ where $R$ is a noetherian local ring, $i$ is a weakly regular (resp. flat complete intersection) local homomorphism and $p$ is surjective. If $B$ is complete, a regular factorization (and so a complete intersection factorization) always exists \cite {AFH}.

We say that a finite module $M\neq0$ over a noetherian local ring $A$ has finite complete intersection dimension in the sense of \cite {AGP} (resp. in the sense of \cite {So}) and use the notation CI-dim$(M)<\infty$ (resp. cidim$(M)<\infty$) if there exists a flat (resp. flat complete intersection) local homomorphism of noetherian local rings $A \rightarrow A'$, and a surjective homomorphism of noetherian local rings $Q \rightarrow A'$ with kernel generated by a regular sequence, such that $pd_Q(M\otimes_AA')<\infty$, where $pd$ denotes projective dimension. For a local homomorphism of noetherian local rings $f:A \rightarrow B$ we say that CI-dim$(f)<\infty$ (resp. cidim$(f)<\infty$) if there exists a regular factorization (resp. complete intersection factorization) $A \rightarrow R \rightarrow \hat{B}$ such that CI-dim$_R(\hat{B})<\infty$ (resp. cidim$_R(\hat{B})<\infty$).

If $fd_A(B)<\infty$, then cidim$(f)<\infty$ (and so also CI-dim$(f)<\infty$). Take a regular factorization $A\rightarrow R \rightarrow \hat{B}$; from the change of rings spectral sequence
$$E^2_{pq} = Tor^R_p(Tor^A_q(R,k),\hat{B}) \Rightarrow Tor^A_{p+q}(k,\hat{B})$$
we obtain $Tor^R_n(R \otimes_Ak,\hat{B}) = Tor^A_n(k,\hat{B}) = 0$ for all $n\gg0$, and then from the spectral sequence
$$E^2_{pq} = Tor^{R \otimes_Ak}_p(Tor^R_q(R \otimes_Ak,\hat{B}),l) \Rightarrow Tor^R_{p+q}(\hat{B},l)$$
we deduce $Tor^R_n(\hat{B},l) = 0$ for all $n\gg0$ since $R \otimes_Ak$ is regular. That is, $pd_R(\hat{B}) = fd_R(\hat{B})<\infty$. This immediately implies cidim$(f)<\infty$.

However, it is clear that cidim$(f)<\infty$ is far from implying $fd_A(B)<\infty$ (take as $A$ a complete intersection local ring and $B$ its residue field).\\

We will need also the following lemma.

\begin{lem}\label{Cohen} {\rm (Essentially \cite [Lemma 1.7]{Av-b})}
Let $(A,k)\xrightarrow{i} (R,l)\rightarrow(D,E)$ be local homomorphisms of noetherian local rings such that $i$ is flat complete intersection. Then the canonical map $H_n(A,D,E) \rightarrow H_n(R,D,E)$ is an isomorphism for all $n\geq 3$ and injective for $n=2$. If $i$ is weakly regular, it is also an isomorphism for $n=2$.
\end{lem}
\begin{proof}
By flat base change $H_n(A,R,E)=H_n(k,R\otimes_Ak,E)$, and by the Jacobi-Zariski exact sequence associated to $k \rightarrow R\otimes_Ak \rightarrow E$ we have $H_n(k,R\otimes_Ak,E)=H_{n+1}(R\otimes_Ak,E,E)=0$ for all $n \geq 2$ by Appendix 9. So the Jacobi-Zariski exact sequence
$$... \rightarrow H_n(A,R,E) \rightarrow H_n(A,D,E) \rightarrow H_n(R,D,E) \rightarrow H_{n-1}(A,R,E) \rightarrow ...$$
gives isomorphisms $H_n(A,D,E)=H_n(R,D,E)$ for all $n\geq 3$ and an exact sequence
$$0 \rightarrow  H_2(A,D,E) \rightarrow H_2(R,D,E) \rightarrow H_1(A,R,E) \xrightarrow{\alpha} H_1(A,D,E) \rightarrow ...$$

The injectivity of $\alpha$ when $i$ is weakly regular follows from the commutative diagram with exact upper row
$$
\xymatrix{  0=H_2(R\otimes_Ak,E,E) \ar[r] & H_1(k,R\otimes_Ak,E) \ar[r] & H_1(k,E,E)  \\
&  H_1(A,R,E) \ar[u]^{\simeq} \ar[r]^{\alpha}  & H_1(A,D,E) \ar[u] }
$$
\end{proof}

The analogues to Avramov's theorems \cite[Theorems 1.4, 1.3]{Av-b} are the following two:

\begin{thm}\label{TheoremA}
Let $f:(A,\mathfrak{m},k) \rightarrow (B,\mathfrak{n},l)$ be a local homomorphism of noetherian local rings such that CI-dim$(f)<\infty$. Let $m \geq 2$ be an integer such that $(m-1)!$ is invertible in $B$. If $H_n(A,B,l)=0$ for some $n$ with $3 \leq n \leq 2m-1$, then $H_n(A,B,l)=0$ for all $n \geq 3$.
\end{thm}

\begin{proof}
By Appendix 8, we may assume that $B$ is complete. Consider a diagram

$$
\xymatrix{ & & Q \ar[dr] & & \\
 & R  \ar[rr] \ar[dr] & & R' \ar[dr] & \\
A \ar[ur] \ar[rr] & & B \ar[rr] & & R'\otimes_RB }
$$
as in the definition of complete intersection dimension. Let $E$ be the common residue field of $R'$ and $R'\otimes_RB$. By Lemma \ref{Cohen} and flat base change, we have isomorphisms
$$H_n(A,B,E) = H_n(R,B,E) = H_n(R',R'\otimes_RB,E)$$
for all $n \geq 2$. Also, from the Jacobi-Zariski exact sequence associated to $Q \rightarrow R' \rightarrow R'\otimes_RB$ and Appendix 9, we obtain
$$H_n(Q,R'\otimes_RB,E) = H_n(R',R'\otimes_RB,E)$$
for all $n \geq 3$.

Therefore, given an integer $n \geq 3$, $H_n(A,B,l) = 0$ if and only if $H_n(Q,R'\otimes_RB,E) = 0$. Since $fd_Q(R'\otimes_RB)<\infty$, by \cite[Theorem 1.4]{Av-b} (by its proof, that is local, or using \cite [4.57]{An1974}) we have $H_n(Q,R'\otimes_RB,E) = 0$ for all $n \geq 2$, and then $H_n(A,B,l) = 0$ for all $n \geq 3$.

\end{proof}

\begin{thm}\label{TheoremB}
Let $f:(A,\mathfrak{m},k) \rightarrow (B,\mathfrak{n},l)$ be a local homomorphism of noetherian local rings such that CI-dim$(f)<\infty$. If $H_n(A,B,l)=0$ for all $n$ sufficiently large, then $H_n(A,B,l)=0$ for all $n \geq 3$.
\end{thm}

\begin{proof}
It follows from \cite[Theorem 4.4]{Av-b} using the same ideas as in the proof of Theorem \ref{TheoremA}.
\end{proof}

\begin{rem}\label{remC}
It should be noted that  working in a similar way with complete intersection dimension, we can also prove part of \cite[Theorem 1.5]{Av-b} in some very particular cases that are not complete intersection. If $f:(A,\mathfrak{m},k) \rightarrow (B,\mathfrak{n},l)$ is a local homomorphism of noetherian local rings with $edim(A)-depth(A) \leq 3$ or $edim(A)-depth(A) = 4$ and $A$ Gorenstein, the ring $R$ in a regular factorization (we can assume that $B$ is complete) $A \rightarrow R \rightarrow B$, inherits the same property. If $H_n(A,B,l)=0$ for all $n$ sufficiently large (and so $H_n(R,B,l)=0$ for all $n$ sufficiently large by Lemma \ref{Cohen}), and we mimic the proof of \cite[Theorem 4.4]{Av-b} with Quillen's spectral sequence \cite[p. 475]{Av-b} instead of the spectral sequence of \cite[Theorem 4.2]{Av-b}, we conclude that the Poincar\'e series
$$\sum_{i \geq 0}{dim_lTor^R_i(B,l)}$$
has radius of convergence $\geq 1$. By \cite{AvAsymp}, this implies for these particular rings that (the virtual projective dimension and so) the complete intersection dimension of the $R$-module $B$ is finite. Therefore reasoning as in the proof of Theorem \ref{TheoremA}, we deduce that $H_n(A,B,l)=0$ for all $n \geq 3$.
\end{rem}

\begin{thm}\label{cidim}
Let
$$
\xymatrix{  (A,\mathfrak{m},k) \ar[r]^u \ar[d]^f & (B,\mathfrak{n},l) \ar[d]^g  \\
(\tilde{A},\mathfrak{\tilde{m}},\tilde{k}) \ar[r]^{\tilde{u}}  & (\tilde{B},\mathfrak{\tilde{n}},\tilde{l}) }
$$
be a commutative square of local homomorphisms of noetherian local rings verifying\\*
(i) $Tor_i^A(\tilde{A},B)=0$ for all $i>0$.\\*
(ii) The homomorphism $H_3(A,B,\tilde{l}) \rightarrow H_3(\tilde{A},\tilde{B},\tilde{l})$ vanishes.\\*
(iii) If $\mathfrak{p}$ is the contraction in $\tilde{A}\otimes_AB$ of the maximal ideal $\mathfrak{\tilde{n}}$ of $\tilde{B}$, then $(\tilde{A}\otimes_AB)_{\mathfrak{p}}$ is a noetherian ring.\\*
(iv) If $\omega$ is the homomorphism $(\tilde{A}\otimes_AB)_{\mathfrak{p}} \rightarrow \tilde{B}$, then cidim$(\omega) <\infty$.\\

Then $H_3(A,B,l)=0$.
\end{thm}

\begin{proof}
We can assume that $\tilde{B}$ is complete. Consider a diagram showing that cidim$(\omega) <\infty$
$$
\xymatrix{ & & Q \ar[dr] & & \\
 & R  \ar[rr] \ar[dr] & & R' \ar[dr] & \\
(\tilde{A}\otimes_AB)_{\mathfrak{p}} \ar[ur] \ar[rr] & & \tilde{B} \ar[rr] & & R'\otimes_R\tilde{B} }
$$

Let $E$ be the residue field of $R'$ and $R'\otimes_R\tilde{B}$. We have
$$H_n((\tilde{A}\otimes_AB)_{\mathfrak{p}},\tilde{B},E) = H_n(R,\tilde{B},E) = H_n(R',R'\otimes_R\tilde{B},E)$$
for all $n \geq 3$, by Lemma \ref{Cohen} and flat base change.

Using \cite [19.21, 20.26, 20.27]{An1974} to interpret \cite{CRAS} in terms of Andr\'e-Quillen homology, since $fd_Q(R'\otimes_R\tilde{B})<\infty$, we have a zero homomorphism in the upper row of the commutative diagram
$$
\xymatrix{H_4(Q,R'\otimes_R\tilde{B},E) \ar[r]^0 \ar[d]^{\simeq} & H_4(Q,E,E) \hphantom{a} \ar @{>->}[r] \ar[d]^{\simeq} & H_4(R'\otimes_R\tilde{B},E,E) \\
H_4(R',R'\otimes_R\tilde{B},E) \ar[r] & H_4(R',E,E) & \\
H_4(R,\tilde{B},E) \ar[u]_{\simeq} & & \\
H_4((\tilde{A}\otimes_AB)_{\mathfrak{p}},\tilde{B},E) \ar[r] \ar[u]_{\simeq} & H_4((\tilde{A}\otimes_AB)_{\mathfrak{p}},E,E) \ar[uu]_{\simeq} & }
$$
where for the isomorphisms we have used again Lemma \ref{Cohen}, and Appendix 4, 9, having in mind that $R \rightarrow R'$ is flat with complete intersection closed fiber. We obtain that the map
$$H_4((\tilde{A}\otimes_AB)_{\mathfrak{p}},\tilde{B},E) \rightarrow H_4((\tilde{A}\otimes_AB)_{\mathfrak{p}},E,E)$$
is zero.

Using hypotheses (i) and (ii) and Appendix 3, 4, we have a commutative diagram
$$
\xymatrix{ H_3(A,B,E) \ar[dr]^0 \ar[rr]^{\simeq} & & H_3(\tilde{A},(\tilde{A}\otimes_AB)_{\mathfrak{p}},E) \ar[dl] \\
& H_3(\tilde{A},\tilde{B},E) & }
$$
showing that the homomorphism $H_3(\tilde{A},(\tilde{A}\otimes_AB)_{\mathfrak{p}},E) \rightarrow H_3(\tilde{A},\tilde{B},E)$ is zero. So from the Jacobi-Zariski exact sequence
$$H_4((\tilde{A}\otimes_AB)_{\mathfrak{p}},\tilde{B},E) \xrightarrow{\delta} H_3(\tilde{A},(\tilde{A}\otimes_AB)_{\mathfrak{p}},E) \xrightarrow{0} H_3(\tilde{A},\tilde{B},E)$$
we deduce that $\delta$ is surjective. Therefore we have a commutative diagram
$$
\xymatrix{ H_4((\tilde{A}\otimes_AB)_{\mathfrak{p}},\tilde{B},E) \ar @{->>}[dr]^{\delta} \ar[rr]^0 & & H_4((\tilde{A}\otimes_AB)_{\mathfrak{p}},E,E) \ar[dl] \\
& H_3(\tilde{A},(\tilde{A}\otimes_AB)_{\mathfrak{p}},E) & }
$$
from which we deduce that $0 = H_3(\tilde{A},(\tilde{A}\otimes_AB)_{\mathfrak{p}},E) = H_3(A,B,E) = H_3(A,B,l)\otimes_lE$.

\end{proof}

We can point out special cases of this theorem as we did in $A-E$ in Section 1. In particular, using Corollary \ref{corollary4} and Theorem \ref{cidim}, we can extend the main results in \cite{BM} from the Frobenius homomorphism to contracting endomorphisms. All works along the same lines, so we content ourselves with the following two criteria for complete intersection homomorphisms.

The reader may consult the definitions and properties of complete intersection homomorphisms and quasi-complete intersection homomorphisms in \cite{Av-b}, \cite{AHS}, where they have been studied in depth. We only need here their characterization in terms of Andr\'e-Quillen homology. If $f:(A,\mathfrak{m},k) \rightarrow (B,\mathfrak{n},l)$ is a local homomorphism of noetherian local rings, $f$ is complete intersection at $\mathfrak{n}$ if and only if $H_n(A,B,l)=0$ for all $n \geq 2$ \cite [Proposition 1.1]{Av-b}, and is quasi-complete intersection at $\mathfrak{n}$ if and only if $H_n(A,B,l)=0$ for all $n \geq 3$ \cite [7.6]{AHS}. We will consider these properties always at the maximal ideals of the rings involved, so again we will suppress ``at $\mathfrak{n}$" from the notation. Note that this terminology is congruent with the notion of flat complete intersection homomorphism defined at the beginning of this section, since if $f$ is flat, then $H_n(A,B,l)=H_n(k,B\otimes_Ak,l)=H_{n+1}(B\otimes_Ak,l,l)$ for $n \geq 2$ by Appendix 4, 6, 7, and $H_{n+1}(B\otimes_Ak,l,l)=0$ for $n \geq 2$ if and only if $B\otimes_Ak$ is a complete intersection ring (Appendix 9).

\begin{cor}
Let $A \xrightarrow{u} B \xrightarrow{v} C$ be local homomorphisms of noetherian local rings such that cidim$(v)<\infty$. Assume that the characteristic of the residue fields is $\neq 2$. If $vu$ is quasi-complete intersection (resp. complete intersection) then $u$ is quasi-complete intersection (resp. complete intersection) and $v$ is quasi-complete intersection.

\end{cor}

\begin{proof}
Let $E$ be the residue field of $C$. We apply Theorem \ref{cidim} to
$$
\xymatrix{ A \ar @{=}[d] \ar[r]^u & B \ar[d]^v \\
A \ar[r]^{vu} & C }
$$
and we obtain $H_3(A,B,E)=0$. Then, from the Jacobi-Zariski exact sequence
$$0=H_4(A,C,E) \rightarrow H_4(B,C,E) \rightarrow H_3(A,B,E)=0$$
we deduce $H_4(B,C,E)=0$, and so from Theorem \ref{TheoremA} we get
$$H_n(B,C,E)=0$$
for all $n \geq 3$, that is, $v$ is quasi-complete intersection.

Again from the same Jacobi-Zariski exact sequence
$$ ... \rightarrow H_{n+1}(B,C,E) \rightarrow H_n(A,B,E) \rightarrow H_n(A,C,E) \rightarrow ...$$
we deduce $H_n(A,B,E)=0$ for all $n \geq 3$ (resp. for all $n \geq 2$).

\end{proof}

In the particular case when $fd_B(C)<\infty$ this result was proved in \cite [5.7.1]{Av-b} for the case of complete intersection and in \cite [7.9]{AHS} for quasi-complete intersection (even in characteristic $2$).\\

\begin{cor}\label{lastcor}
Let $u:(A,\mathfrak{m},k) \rightarrow (B,\mathfrak{n},l)$ be a local flat homomorphism of noetherian local rings and $\omega :B\otimes_AB \rightarrow B$ the multiplication map. Let $\mathfrak{p} = \omega^{-1}(\mathfrak{n})$ and assume that $(B\otimes_AB)_{\mathfrak{p}}$ is noetherian. If cidim$_{(B\otimes_AB)_{\mathfrak{p}}}(B)<\infty$, then $u$ is complete intersection (at $\mathfrak{n}$).

\end{cor}
\begin{proof}
By Theorem \ref{cidim} applied to
$$
\xymatrix{A \ar[d]^u \ar[r]^u & B \ar @{=}[d] \\
B \ar @{=}[r] & B }
$$
we obtain $H_3(A,B,l)=0$, and since $u$ is flat, we deduce from \cite [17.2]{An1974} that $H_n(A,B,l)=0$ for all $n \geq 2$.

\end{proof}

\begin{rem}
Similarly to Remark \ref{augmented}, this result is also valid more generally for a noetherian supplemented $B$-algebra $S$ instead of $B\otimes_AB$. If cidim$_{S_{\mathfrak{p}}}(B)<\infty$, where $\mathfrak{p}$ is the contraction in $S$ of the maximal ideal of $B$, we apply Theorem \ref{cidim} to the square
$$
\xymatrix{B \ar @{=}[d] \ar[r] & S_{\mathfrak{p}} \ar[d] \\
B \ar @{=}[r] & B }
$$
We obtain $H_3(B,S_{\mathfrak{p}},l)=0$, and so from the Jacobi-Zariski exact sequence associated to $B \rightarrow S_{\mathfrak{p}} \rightarrow B$ we deduce $H_4(S_{\mathfrak{p}},B,l)=0$. By Theorem \ref{TheoremA} (in characteristic $\neq 2$), $H_n(S_{\mathfrak{p}},B,l)=0$ for all $n \geq 3$, and then from \cite [4.57]{An1974} we deduce $H_n(S_{\mathfrak{p}},B,-)=0$ for all $n \geq 3$.

In fact, the hypothesis on the characteristic is not necessary since by \cite [Theorem I]{AI}, for a noetherian supplemented $B$-algebra $S$, $H_4(S,B,-)=0$ implies $H_n(S,B,-)=0$ for all $n \geq 3$ (a direct proof of this fact was communicated to me by Rodicio: we can assume that the rings are local and complete; take a regular factorization $B \rightarrow R\rightarrow S$; we have $H_3(B,S,-)=H_4(S,B,-)=0$ and so $H_3(R,S,-)=0$ by Lemma \ref{Cohen}; using this same Lemma, we have $H_n(R,B,-)=H_n(B,B,-)=0$ for all $n \geq 2$; in the commutative diagram with exact row
$$
\xymatrix{H_3(S,B,l) \ar[d] \ar[r] & H_2(R,S,l) \ar[r] \ar @{=}[d] & H_2(R,B,l)=0 \\
H_3(S,l,l) \ar[r]^{\alpha} & H_2(R,S,l) }
$$
we have $\alpha=0$ by \cite [Corollary 6]{Ro3} and then $H_2(R,S,l)=0$; this implies $H_n(R,S,l)=0$ for all $n \geq 2$ by Appendix 9 and then $H_n(B,S,l)=0$ for all $n \geq 2$ by Lemma \ref{Cohen}; since $H_{n+1}(S,B,l)= H_n(B,S,l)$ we deduce $H_n(S,B,-)=0$ for all $n \geq 3$).

\end{rem}

\appendix
\section{Some results on Andr\'e-Quillen homology used in the proofs}

Associated to a homomorphism of (commutative) rings $f:A \rightarrow B$ and to a $B$-module $M$ we have (Andr\'e-Quillen) homology $B$-modules $H_n(A,B,M)$ for all integers $n \geq 0$, which are functorial in all three variables and satisfy the following properties:\\

1. If $B=A/I$, then $H_0(A,B,M)=0$  \cite [4.60]{An1974}. \\

2. If $0 \rightarrow M' \rightarrow M \rightarrow M'' \rightarrow 0$ is an exact sequence of $B$-modules, we have a natural exact sequence \cite [3.22]{An1974}
\begin{align*}
... \rightarrow H_{n+1}(A,B,M'')\to \\
H_n(A,B,M') \rightarrow H_n(A,B,M)\rightarrow H_n(A,B,M'')\rightarrow \\
H_{n-1}(A,B,M') \rightarrow \hspace{5mm}
... \hspace{5mm} \rightarrow H_0(A,B,M'')\rightarrow 0\\
\end{align*}

3. (Localization) Let $u: A \rightarrow B$ be a ring homomorphism, $T$ a multiplicative subset of $B$, $S$ a multiplicative subset of $A$ such that $u(S)\subset T$, and $M$ a $B$-module. Then \cite [4,59, 5.27]{An1974}
\begin{flushleft}
$T^{-1}H_n(A,B,M) = H_n(A,B,T^{-1}M) = $
\end{flushleft}
\begin{flushright}
$H_n(A,T^{-1}B,T^{-1}M) = H_n(S^{-1}A,T^{-1}B,T^{-1}M).$
\end{flushright}

\hspace{1cm}

4. (Base change) Let $A \rightarrow B$, $A \rightarrow C$ be ring homomorphisms such that $Tor_i^A(B,C)=0$ for all $i>0$, and let $M$ be a $B\otimes_AC$-module. Then $H_n(A,B,M)=H_n(C,B\otimes_AC,M)$ for all $n$ \cite [4.54]{An1974}. \\

5. Let $B$ be an $A$-algebra, $C$ a $B$-algebra and $M$ a flat $C$-module. Then $H_n(A,B,M)=H_n(A,B,C) \otimes_CM$ for all $n$ \cite [3.20]{An1974}. \\

6. (Jacobi-Zariski exact sequence) If $A \rightarrow B \rightarrow C$ are ring homomorphisms and $M$ is a $C$-module, we have a natural exact sequence \cite [5.1]{An1974}

\begin{align*}
... \rightarrow H_{n+1}(B,C,M)\to \\
H_n(A,B,M) \rightarrow H_n(A,C,M)\rightarrow H_n(B,C,M)\rightarrow \\
H_{n-1}(A,B,M) \rightarrow \hspace{5mm}
... \hspace{5mm} \rightarrow H_0(B,C,M)\rightarrow 0\\
\end{align*}

7. If $K\rightarrow L$ is a field extension and $M$ an $L$-module, we have $H_n(K,L,M)=0$ for all $n \geq 2$ \cite [7.4]{An1974}. So if $A\rightarrow K \rightarrow L$ are ring homomorphisms with $K$ and $L$ fields, from 6 we obtain $H_n(A,K,L) = H_n(A,L,L)$ for all $n \geq 2$, which, using 5, gives $H_n(A,K,K) \otimes_KL=H_n(A,L,L)$ for all $n \geq 2$. \\

8. If $(A,\mathfrak{m},k)$ is a noetherian local ring and $\hat{A}$ is its $\mathfrak{m}$-completion, then $H_n(A,k,k)=H_n(\hat{A},k,k)$ for all $n \geq 0$ \cite [10.18]{An1974}. As a consequence, if $(A,\mathfrak{m},k) \rightarrow (B,\mathfrak{n},l)$ is a local homomorphism of noetherian local rings, then taking the map between the Jacobi-Zariski exact sequences associated to $A \rightarrow B \rightarrow l$ and $\hat{A} \rightarrow \hat{B} \rightarrow l$, we deduce that $H_n(A,B,l)=H_n(\hat{A},\hat{B},l)$.\\

9. If $I$ is an ideal of a noetherian local ring $(A,\mathfrak{m},k)$, then the following are equivalent: \\*
(i) $I$ is generated by a regular sequence \\*
(ii) $H_2(A,A/I,k)=0$ \\*
(iii) $H_n(A,A/I,M)=0$ for any $A/I$-module $M$ for  all $n \geq 2$ \cite [6.25]{An1974}. \\

In particular, a noetherian local ring $(A,\mathfrak{m},k)$ is regular if and only if $H_2(A,k,k)=0$.

Similarly, a noetherian local ring $(A,\mathfrak{m},k)$ is complete intersection if and only if $H_3(A,k,k)=0$ if and only if $H_n(A,k,k)=0$ for  all $n \geq 3$ \cite [6.27 and its proof]{An1974}.\\

Finally, two important results by Andr\'e and Avramov respectively that play a key role in this paper:\\

10. A local homomorphism of noetherian local rings $(A,\mathfrak{m},k)\rightarrow(B,\mathfrak{n},l)$ is formally smooth (for the $\mathfrak{n}$-adic topology) in the sense of \cite [0$_{IV}$19.3.1]{EGAIV1} if and only if $H_1(A,B,l)=0$ (essentially \cite [16.17]{An1974}; use \cite [3.21]{An1974}). \\

11. The main result in \cite {CRAS} can be read in terms of Andr\'e-Quillen homology as follows (see \cite [19.21, 20.26]{An1974}): if $(A,\mathfrak{m},k)\rightarrow(B,\mathfrak{n},l)$ is a local homomorphism of finite flat dimension of noetherian local rings, then the homomorphism $H_2(A,l,l) \rightarrow H_2(B,l,l)$ is injective. \\


\end{document}